# EFFICIENT RARE-EVENT SIMULATION FOR THE MAXIMUM OF HEAVY-TAILED RANDOM WALKS[1]

BY JOSE BLANCHET AND PETER GLYNN

*Harvard University and Stanford University*

Let $(X_n : n \geq 0)$ be a sequence of i.i.d. r.v.'s with negative mean. Set $S_0 = 0$ and define $S_n = X_1 + \cdots + X_n$. We propose an importance sampling algorithm to estimate the tail of $M = \max\{S_n : n \geq 0\}$ that is *strongly efficient* for both light and heavy-tailed increment distributions. Moreover, in the case of heavy-tailed increments and under additional technical assumptions, our estimator can be shown to have *asymptotically vanishing relative variance* in the sense that its coefficient of variation vanishes as the tail parameter increases. A key feature of our algorithm is that it is state-dependent. In the presence of light tails, our procedure leads to Siegmund's (1979) algorithm. The rigorous analysis of efficiency requires new Lyapunov-type inequalities that can be useful in the study of more general importance sampling algorithms.

**1. Introduction.** In this paper we consider the problem of efficient simulation of first-passage time probabilities for heavy-tailed random walks (r.w.'s). More precisely, suppose that $(S_n : n \geq 0)$ is the r.w. generated by the sequence of independent and identically distributed (i.i.d.) random variables (r.v.'s) $X = (X_n : n \geq 1)$ (i.e., $S_n = S_{n-1} + X_n$ with $S_0 = 0$). We assume that $EX_n < 0$. Define $M = \max\{S_n : n \geq 0\}$ and $\tau(b) = \inf\{n \geq 0 : S_n > b\}$. We are interested in developing efficient simulation methodology to compute

$$P(\tau(b) < \infty) = P(M > b), \tag{1}$$

when $b$ is large (i.e., the event $\{M > b\}$ is rare) and $X_1$ is heavy-tailed.

Received January 2007; revised September 2007.
[1]Supported in part by NSF Grant DMS-05-95595.
*AMS 2000 subject classifications.* Primary 60G50, 60J05, 68W40; secondary 60G70, 60J20.
*Key words and phrases.* State-dependent importance sampling, rare-event simulation, heavy-tails, Lyapunov bounds, random walks, single-server queue, change-of-measure.

We say that an unbiased simulation estimator $R(b)$ for $P(M > b)$ is *strongly efficient* if

$$\sup_{b>0} ER(b)^2 / P(M > b)^2 < \infty.$$

Strong efficiency implies that the number of simulation runs required to estimate $P(M > b)$ to a given relative accuracy is bounded in $b$. A weaker criterion is *logarithmic efficiency*, which implies that the number of replications required to estimate $P(M > b)$ with a given relative accuracy grows at rate $o(|\log P(M > b)|)$; see Asmussen and Glynn (2007), Juneja and Shahabuddin (2006) or Bucklew (2004), Section 5.2, for a discussion of efficiency in rare-event simulation. A strongly efficient estimator is said to exhibit *asymptotically vanishing relative error* when $ER(b)^2 \sim P(M > b)^2$ as $b \nearrow \infty$ (or, equivalently, when the coefficient of variation vanishes as $b \nearrow \infty$).

In this paper we develop an implementable state-dependent importance sampling algorithm that can be rigorously proved to possess asymptotically vanishing relative error. By "state-dependent," we mean that the importance sampling algorithm generates the next increment of the random walk from a distribution that depends on the walk's current state (i.e., location). This is the first strongly efficient algorithm that has been developed for estimating the tail of $M$ in the presence of general heavy-tailed increment distributions. Prior efficient algorithms require the increment distribution to be of $M/G/1$ type with regularly varying or Weibull type right tails.

A key idea is that our importance distribution is state-dependent. There is a long history of applications of state-dependent importance sampling to simulation problems. Perhaps the first related contributions are those by Hammersley and Morton (1954) and Rosenbluth and Rosenbluth (1955) in the context of molecular simulation; see also the text by Liu (2001) for applications of sequential importance sampling in various scientific contexts. However, a general framework for rigorous analysis of these types of algorithms is still under development. In a sequence of recent papers, Paul Dupuis and Hui Wang [see, e.g., Dupuis and Wang (2004)] have proposed a general methodology that can be applied in the presence of large deviations theory for light-tailed systems. Our paper contributes to this general literature by developing Lyapunov-type inequalities (see Theorem 2) that are particularly useful for the analysis of state-dependent algorithms.

The general theory of importance sampling establishes that the theoretically optimal importance distribution (having zero variance) involves sampling from the conditional distribution of the random walk given $\{\tau(b) < \infty\}$. Under this conditional distribution, the random walk has increment distributions that are state-dependent. However, we cannot implement this zero variance sampling scheme because the state-dependent increment distribution requires explicit knowledge of the function $u^*(\cdot) = P(\tau(\cdot) < \infty)$.

Our approach involves using asymptotic approximations for $u^*(\cdot)$ to obtain an implementable state-dependent change-of-measure that closely approximates the true conditional distribution. In the current $G/G/1$ setting, the asymptotic approximation for $u^*(\cdot)$ is

$$u^*(b) = P(M > b) \sim \frac{1}{|EX|} \int_b^\infty P(X > s)\,ds, \tag{2}$$

as $b \nearrow \infty$. An important step in our approach is to use (2) in order to construct a function $v(\cdot)$ such that

$$Ev(b + X) - v(b) = o(P(X > b)) \tag{3}$$

as $b \nearrow \infty$. Note that if $v = u^*$, the above difference vanishes. The above convergence rate [namely, that associated with (3)] is a convenience in developing our simulation algorithm, but is not necessary (see Proposition 3 and Theorem 3). We show that a $v(\cdot)$ satisfying (3) can be constructed using (2) whenever $X$ belongs to the class $S^*$ of heavy-tailed distributions—which is slightly smaller than the class of subexponential distributions but includes regularly varying, Weibull, lognormal and many more distributions as special cases; see Assumption A in Section 3 for a precise definition.

The problem that we address here is motivated by applications in queueing and insurance. The distribution of $M$ is of great interest in queueing theory as it coincides with the steady-state waiting time distribution of the single-server $G/G/1$ queue. In addition, the first passage time probability displayed in (1) is of central interest in the context of insurance risk. In particular, such a first passage time probability can be interpreted as the probability that an insurer receiving premiums at a constant rate is eventually ruined when subject to a renewal arrival process of i.i.d. claims. When the claim distribution is heavy-tailed, the resulting calculation is exactly of the type discussed in this paper. Statistical evidence suggests that such heavy-tailed distributions frequently arise in practice and are a convenient vehicle for capturing many of the key stylized features that are present in observed claim sizes [see, e.g., Embrechts, Klüppelberg and Mikosch (1997) and Adler, Feldman and Taqqu (1998)].

The first efficient rare event simulation algorithm for the tail of $M$ was suggested by Siegmund (1976), who was motivated by the first passage time interpretation displayed in (1) and its connection to one-sided sequential probability ratio tests in the context of statistical sequential analysis. Siegmund's algorithm applies only to light-tailed r.w.'s and involves an importance distribution corresponding to a r.w. with state-independent increments. Our proposed strongly efficient algorithm is consistent with recent results of Bassamboo, Juneja and Zeevi (2006), who show that no state independent efficient importance sampling algorithm for computing (1) can exist in the (regularly varying) heavy-tailed setting. Another key feature that is present in the light-tailed context is the ability to fully leverage the existing theory of large deviations. A complicating factor in the heavy-tailed setting is that the large deviations literature is not applicable to such problems. Asmussen, Binswanger and Hojgaard (2000) provide a number of examples and counterexamples to illustrate the additional difficulties that arise in the heavy-tailed environment.

As noted above, rare-event simulation algorithms for heavy-tailed distributions have been previously developed in the context of the $M/G/1$ queue. The first logarithmically efficient simulation algorithm for estimation of (1) was given in Asmussen and Binswanger (1997) and was based on the idea of conditional Monte Carlo (and not importance sampling). Logarithmic efficiency for their algorithm was established for regularly varying tails and was shown to fail for Weibull-type heavy tails. Subsequently, Asmussen, Binswanger and Hojgaard (2000) developed simulation estimators for the $M/G/1$ queue based on importance sampling ideas that are provably logarithmically efficient for both regularly varying and Weibull-type tails. Juneja and Shahabuddin (2002) also developed logarithmically efficient importance sampling schemes based on a suitable twisting of the $M/G/1$ service time distribution's hazard rate. More recently, Asmussen and Kroese (2006) proposed other logarithmically efficient importance sampling algorithms for the $M/G/1$ queue that seem to have excellent performance in practice. In addition, they developed a conditional Monte Carlo estimator that is strongly efficient for both regularly varying tails and certain Weibull type heavy-tails. Dupuis, Leder and Wang (2006) proposed a state-dependent importance sampling algorithm that is strongly efficient for a regularly varying $M/G/1$ queue. All the above algorithms take advantage of the fact that the ladder height distribution for the $M/G/1$ queue is explicitly known. In contrast, no such explicit computations are possible for the class of $G/G/1$ models considered here. This significantly complicates both the development and the theoretical analysis of efficient rare-event algorithms for this class of problems. Indeed, we developed our Lyapunov bounds largely in order to provide a suitable verification tool for bounding the variances (as required to establish strong efficiency) for the algorithm considered here. More recently, Blanchet, Glynn and Liu (2007) have used this Lyapunov technique to study an alternative importance sampling algorithm for the $G/G/1$ queue that is based on mixture sampling (rather than on a Markovian importance sampler having a transition kernel based explicitly on the approximation $v$ as is the case here). This alternative algorithm, while typically simpler to implement than the approach described here (because generating increments from a mixture distribution is usually easier than the variate generation schemes required here), applies only to regularly varying distributions (rather than the class $S^*$ covered by this paper's algorithm).

The paper is organized as follows. Section 2 introduces a general technique to study efficient state-dependent importance sampling algorithms for computing first passage time probabilities of general state-space Markov chains and recovers Siegmund's algorithm as a direct application of the basic ideas underlying our procedure. Section 3 introduces the precise technical assumptions under which we develop our methodology and provides both the proof of strong efficiency for our importance sampling estimator and establishes its asymptotically vanishing relative error property. In Section 4 we discuss computational complexity issues associated with our algorithm, leading us to a study of the number of variate generations required to terminate our procedure. Additional practical observations and some results on simulation experiments are given in our final section.

**2. Efficient importance samplers for exit probabilities.** The problem of computing the level crossing probability (1) can be viewed as a special case of computing an exit probability. To be specific, let $Y = (Y_n : n \geq 0)$ be a $\mathcal{X}$-valued Markov chain (with stationary transition probabilities) and let $P_y(\cdot)$ and $E_y(\cdot)$ be the probability distribution and expectation operator on the path-space of $Y$, conditional on $Y_0 = y$. For $B \subseteq \mathcal{X}$, let $T = \inf\{n \geq 0 : Y_n \in B\}$ be the exit time from $B^c$. For $A \subseteq B$, the probability $u^*(y) = P_y(Y_T \in A, T < \infty)$ is called an "exit probability" (all the sets considered here are assumed measurable). Note that the level crossing probability (1) is the special case in which $Y$ is given by the r.w. $(S_n : n \geq 0)$, $\mathcal{X} = [-\infty, \infty)$, $B = \{-\infty\} \cup (b, \infty)$, $A = (b, \infty)$ and $y = 0$. Because of the translation invariance of r.w., studying this problem as $b \nearrow \infty$ is equivalent to fixing $B = \{-\infty\} \cup [0, \infty)$, $A = (0, \infty)$, setting $y = -b$ and letting $b \nearrow \infty$. With $B$ and $A$ fixed in this way, our goal is to efficiently compute $u^*(-b)$ as $b \nearrow \infty$. This reformulation of the problem will form the basis of our analysis in the remainder of the paper.

The following result is easily proved [see, e.g., Meyn and Tweedie (1993)].

PROPOSITION 1. *The function $u^* = (u^*(y) : y \in B^c)$ is the minimal nonnegative solution to*

$$u(y) = \int_{\mathcal{X}} P_y(Y_1 \in dz) u(z), \qquad y \in B^c,$$

*subject to the boundary conditions that $u(z) = 1$ for $z \in A$ and $u(z) = 0$ for $z \in B \cap A^c$.*

As mentioned in the Introduction, the zero-variance importance distribution for computing $u^*(y)$ is that associated with the conditional distribution $P_y(\cdot | Y_T \in A, T < \infty)$. Let $\mathcal{F}_n = \sigma(Y_0, \ldots, Y_n)$ for $n \geq 0$. Our next result characterizes this conditional distribution.

THEOREM 1. *Suppose that $u^*(y) > 0$ for $y \in B^c$. Then, for each nonnegative $\mathcal{F}_T$-measurable r.v. $\Lambda$,*

$$E_y[\Lambda | Y_T \in A, T < \infty] = E_y^* \Lambda,$$

*where $E_y^*(\cdot)$ is the expectation operator under which $Y$ is a Markov chain having one-step transition kernel*

$$P^*(y, dz) = P_y(Y_1 \in dz) \frac{u^*(z)}{u^*(y)},$$

*for $y \in B^c$, $z \in \mathcal{X}$.*

PROOF. Note that $I(T = n)\Lambda = \lambda_n(Y_0, \ldots, Y_n)$ for some (measurable) function $\lambda_n : \mathcal{X}^{n+1} \to [0, \infty)$. Therefore,

$$\frac{E_y[\Lambda; T = n, Y_T \in A, T < \infty]}{u^*(y)}$$
$$= \int_{B^c \times \cdots \times B^c \times A} \frac{\lambda_n(y, z_1, \ldots, z_n) u^*(z_n) P(y, dz_1) \cdots P(z_{n-1}, dz_n)}{u^*(y)}$$
$$= \int_{B^c \times \cdots \times B^c \times A} \lambda_n(y, z_1, \ldots, z_n) \frac{P(y, dz_1) u^*(z_1)}{u^*(y)} \frac{P(z_1, dz_2) u^*(z_2)}{u^*(z_1)}$$
$$\times \cdots \times \frac{P(z_{n-2}, dz_{n-1}) u^*(z_{n-1})}{u^*(z_{n-2})} \frac{P(z_{n-1}, dz_n) u^*(z_n)}{u^*(z_{n-1})}$$
$$= E_y^*[\Lambda; T = n].$$

Summing over $n$, we conclude that

$$E[\Lambda | Y_T \in A, T < \infty] = E^*[\Lambda; T < \infty].$$

Letting $\Lambda = 1$ establishes that $P_y^*(T < \infty) = 1$, proving the result. □

This theorem makes clear that the zero-variance importance sampling distribution for computing (1) corresponds to a random walk in which the increments have a state-dependent distribution. The above result suggests that a good importance sampling distribution can be obtained by simulating $Y$ under transition dynamics that closely approximate those induced by the zero-variance importance distribution's transition kernel $P^*$.

Suppose that $Q$ is the Markov transition kernel chosen by the simulator to compute the exit probability $u^*(y) = P_y(Y_T \in A, T < \infty)$ via importance sampling. Assume that $(Q(y, dz) : y, z \in B^c \cup A)$ can be represented as

$$Q(y, dz) = r(y, z)^{-1} P_y(Y_1 \in dz) I(y \in B^c, z \in B^c \cup A)$$
$$+ \delta_y(dz) I(y \in A, z \in A)$$

for some positive function $r(\cdot)$. Note that

$$P_y(Y_T \in A, T = n)$$
$$= E_y^Q\left[I(Y_T \in A, T = n)\prod_{j=1}^{T} r(Y_{j-1}, Y_j)\right],$$

where $E_y^Q(\cdot)$ is the expectation operator under which $Y$ evolves according to the transition kernel $Q$, conditional on $Y_0 = y$. Summing over $n$, we conclude that $u^*$ can be represented as

$$u^*(y) = E_y^Q\left[I(T < \infty)\prod_{j=1}^{T} r(Y_{j-1}, Y_j)\right].$$

An important step in any theoretical analysis of the estimator

$$(4) \qquad R = I(T < \infty)\prod_{j=1}^{T} r(Y_{j-1}, Y_j)$$

is to bound its variance. The variance, conditional on $Y_0 = y$, is given by $s^*(y) - u^*(y)^2$, where $s^*(y) = E_y^Q R^2$. Since only $s^*(\cdot)$ depends on the choice of the importance distribution, we focus on bounding this quantity.

THEOREM 2.

(i) *The function $s^* = (s^*(y) : y \in B^c)$ is the minimal nonnegative solution to*

$$s(y) = \eta(y) + \int_{B^c} K(y, dz)s(z),$$

*for $y \in B^c$, where*

$$\eta(y) = \int_A r(y, z)P_y(Y_1 \in dz),$$
$$K(y, dz) = r(y, z)P_y(Y_1 \in dz),$$

*for $y, z \in B^c$.*

(ii) *The function $s^*$ is given by*

$$s^* = \sum_{n=0}^{\infty} K^n \eta,$$

*where $K^n(y, dz) = \int_{B^c} K^{n-1}(y, dy_1)K(y_1, dz)$ for $n \geq 1$, $K^0(y, dz) = \delta_y(dz)$ and $(K^n \eta)(y) = \int_{B^c} K^n(y, dz)\eta(z)$.*

(iii) *Suppose that $h = (h(y): y \in B^c)$ is a finite-valued nonnegative function for which*

$$(Kh)(y) \leq h(y) - \eta(y) \tag{5}$$

*for $y \in B^c$. Then, $s^*(y) \leq h(y)$ for $y \in B^c$.*

PROOF. Part (ii) follows by expanding $E_y^Q[R^2 I(T = n)]$ and summing over $n$ using Fubini's theorem. Part (i) follows easily from (ii).

For part (iii), first note that $Kh$ must be finite-valued by virtue of (5). Induction based on applying $K^n$ to both sides of (5) establishes that $K^n h$ is finite-valued for $n \geq 1$. By applying $K^n h$ to (5) and using the fact that $K^n h$ is finite-valued for $n \geq 1$, we conclude that $K^n \eta \leq K^n h - K^{n+1} h$ for $n \geq 0$. Summing over $0 \leq n \leq m$ and using the nonnegativity of $h$, we obtain the bound

$$\sum_{n=0}^{m} K^n \eta \leq h - K^{m+1} h \leq h.$$

The result follows by sending $n \nearrow \infty$ and using part (iii). □

We call the function $h(\cdot)$ a Lyapunov function and refer to bounds based on part (iii) of Theorem 2 as Lyapunov bounds on the second moment.

Returning to the exit probability computations, suppose that $v = (v(y): y \in \mathcal{X})$ is chosen by the simulator to be a good approximation to $u^* = (u^*(y): y \in \mathcal{X})$. In view of Theorem 2 above, it is then natural to consider simulating $Y$ via the transition kernel

$$Q(y, dz) = P(y, dz) \frac{v(z)}{w(y)} \tag{6}$$

(for $y \in B^c$, $z \in B^c \cup A$), where $w(y)$ is the normalization constant given by

$$w(y) = \int_{B^c \cup A} P(y, dz) v(z)$$

(assumed to be finite). In this case, $r(y, z) = w(y)/v(z)$. The following result provides a Lyapunov bound on the second moment $s^*(\cdot)$ that is specifically suited to this setting.

PROPOSITION 2. *Assume that $w(y) > 0$ for $y \in B^c$ and suppose that there exists a finite-valued function $h: B^c \cup A \longrightarrow [\varepsilon, \infty)$ satisfying*

$$w(y) \int v(z) h(z) P(y, dz) \leq h(y) v(y)^2, \tag{7}$$

*for $y \in B^c$. If $h(z) \geq 1$ for $z \in A$ and $v(z) \geq \kappa > 0$ for $z \in A$, then $s^*(y) \leq \varepsilon^{-1} \kappa^{-2} v(y)^2 h(y)$.*

PROOF. Put $\widetilde{h}(\cdot) = \kappa^{-2} h(\cdot) v^2(\cdot)$ and note that (7) is equivalent to assuming that

$$(8) \qquad (K\widetilde{h})(y) \leq \widetilde{h}(y) - \kappa^{-2} w(y) \int_A P(y, dz) w(y) h(z)$$

for $y \in B^c$. But

$$\eta(y) = \int_A P(y, dz) \frac{w(y)}{v(z)} \leq \int_A \kappa^{-2} P(y, dz) w(y) v(z)$$
$$\leq \kappa^{-2} w(y) \varepsilon^{-1} \int_A P(y, dz) v(z) h(z),$$

so that (8) implies that

$$(K\widetilde{h})(y) \leq \widetilde{h}(y) - \eta(y)$$

for $y \in B^c$. We now apply part (iii) of Theorem 2 to complete the proof. □

Suppose that $v(\cdot)$ has been chosen by the simulator to be within a constant multiple of $u^*(\cdot)$, as occurs whenever $v(\cdot)$ has the same asymptotic behavior as $u^*(\cdot)$. In this case, it follows that the importance sampling algorithm based on $r(y, z) = w(y)/v(z)$ has bounded relative variance [i.e., the ratio of the variance to the square of $u^*(x)$] across $B^c$ whenever the function $h$ of Proposition 2 can be chosen to be bounded. On the other hand, if $h$ grows at a suitable rate [e.g., $h(y) = |\log(v(y))|^{1/2}$], the logarithmic efficiency of the importance sampler can be assured.

To illustrate, consider the problem of estimating

$$u^*(-b) = P(\tau(0) < \infty | S_0 = -b)$$

for $b > 0$ in the light-tailed setting. In particular, suppose that there exists a positive root $\theta^*$ of $E[\exp(\theta^* X_1)] = 1$ for which $E[X_1 \exp(\theta^* X_1)] < \infty$. If $X_1$ is nonlattice, then it is known that

$$u^*(-b) \sim c \exp(-\theta^* b)$$

for some positive constant $c$; see, for example, Asmussen (2003), page 365. The natural choice for $v$ is, of course, $v(z) = \exp(\theta^* z)$, in which case $w(z) = \exp(\theta^* z)$. If we put $h(y) = 1$ for $y \in \mathbb{R}$, Proposition 2 applies, yielding the bound

$$s^*(-b) \leq \exp(-2\theta^* b).$$

Hence, this importance sampling algorithm [which is precisely the one proposed by Siegmund (1976)] is strongly efficient.

**3. Elements of our algorithm for heavy-tailed r.w.'s.** We shall explore how to adapt the ideas discussed in the previous sections to the case of a random walk with heavy-tailed increment distributions. We need the following definitions. Set $X^+ = \max(X, 0)$ and $X^- = \max(-X, 0)$.

DEFINITION 1. A nonnegative r.v. $Z$ is said to be *subexponential* if

$$P(Z_1 + Z_2 > t) \sim 2P(Z > t),$$

as $t \nearrow \infty$ where $Z_1$ and $Z_2$ are independent copies of $Z$. A r.v. $X$ is said to be subexponential if $X^+$ is subexponential.

DEFINITION 2. A nonnegative r.v. $Z$ belongs to the family $S^*$ if

$$2EZP(Z > t) \sim \int_0^t P(Z > t - s)P(Z > s)\,ds$$

as $t \nearrow \infty$. In addition, a r.v. $X$ is in $S^*$ if $X^+$ is in $S^*$.

DEFINITION 3. A r.v. $X$ is said to possess a *long tail* if for every constant $a \in \mathbb{R}$

$$P(X > t + a) \sim P(X > t)$$

as $t \nearrow \infty$.

It can be shown that if $Z$ is in $S^*$, then it must be subexponential. Also, any subexponential r.v. possesses a long tail. The class $S^*$ of random variables includes, as particular cases, regularly varying, lognormal and Weibull-type distributions among many others. For more on the specific properties of various types of heavy-tailed distributions, see Embrechts, Klüppelberg and Mikosch (1997), Section 1.4.

The following assumption will be imposed throughout the rest of the paper.

ASSUMPTION A. Assume that $X_n^+$ belongs to $S^*$, that is,

$$2EX_n^+ P(X_n > t) \sim \int_0^t P(X_n > t - s)P(X_n > s)\,ds$$

as $t \nearrow \infty$.

If $X$ belongs to $S^*$, then both the distribution of $X$ and its integrated tail

$$\int_x^\infty \frac{P(X > s)}{EX^+}\,ds$$

are subexponential [see Asmussen (2003), Section 10.9]. Under Assumption A, it is known [see, e.g., Asmussen (2003), page 296] that

$$(9) \qquad u^*(-b) = P(\tau(0) < \infty | S_0 = -b) \sim \frac{-1}{EX} \int_b^\infty P(X > t)\, dt$$

as $b \nearrow \infty$. The previous result is also known in the literature as the Pakes–Veraberbeke theorem.

The natural strategy is to use this approximation to construct an appropriate importance sampling transition kernel $Q(x, dy)$ [defined in (6)] by means of a function $v(\cdot)$ that mimics the behavior of $u^*(\cdot)$. An important estimate in the efficiency analysis of our importance sampling scheme involves the behavior of $v(y) - w(y)$ as $y \searrow -\infty$, where $w(y) = Ev(y + X)$. As we indicated earlier, if one selects $v = u^*$, then the difference $v(y) - w(y)$ vanishes. Thus, it is natural to expect that the asymptotic behavior of this difference will play an important role in the performance of the importance sampling estimator. As we shall see, in order to guarantee strong efficiency of the importance sampling estimator, it suffices to select $v(\cdot)$ so that $v(y) - w(y) = o(P(X > -y))$ as $y \searrow -\infty$.

Recent estimates by Borovkov and Borovkov (2001) under regularly varying or semiexponential assumptions provide asymptotics to $u^*(y)$ that hold with an error of order $o(P(X > -y))$ as $y \searrow -\infty$. Under these assumptions, Borovkov and Borovkov (2001) add an additional term to (9) of order $O(P(X > -y))$ to the approximation (9) which yields an error rate $o(P(X > -y))$ as $y \searrow -\infty$.

Given the form of (9), it may be surprising at first sight that making use only of approximation (9) and assuming only that the distribution of $X$ belongs to the class $S^*$ one can easily construct $v(\cdot)$ that actually achieves an error of order $o(P(X > -y))$ for the difference $v(y) - w(y)$ as $y \searrow -\infty$. In fact, as we shall prove in our next proposition, $v(-t)$ can be defined as the tail probability of a nonnegative random variable $Z$ such that

$$(10) \qquad P(Z > t) = \min\left[-(EX)^{-1} \int_t^\infty P(X > s)\, ds, 1\right]$$

for $t > 0$ [this may imply $P(Z = 0) > 0$]. Then, we write $v(y) = P(Z > -y)$ for all $y \in \mathbb{R}$. Note that if we could pick $u^* = v$, this would correspond to choosing $Z = M$. Given our representation for $v(\cdot)$ as a tail probability, we can write

$$w(y) = E[v(y + X)] = P(X + Z > -y).$$

The next result shows that this choice of $v(\cdot)$ has the indicated convergence rate for the difference $v(y) - w(y)$. However, for the purpose of our efficiency analysis, it is the second part of the following result, namely, inequality (11), which we shall invoke.

PROPOSITION 3. *Under Assumption A,*
$$w(y) - v(y) = o(P(X > -y))$$
*as $y \searrow -\infty$. Consequently, for each $\gamma \in (0,1)$, there exists $a_*(\gamma) \in (-\infty, 0]$ such that, for all $y \leq a_*(\gamma)$,*

(11) $$-\gamma \leq \frac{v(y)^2 - w(y)^2}{P(X > -y)w(y)}.$$

PROOF. We must show that
$$P(X + Z > t) - P(Z > t) = o(P(X > t))$$
as $t \nearrow \infty$. Note that
$$P(X + Z > t) = P(X + Z > t, Z > t) + P(X + Z > t, Z \leq t)$$
$$= P(Z > t) - P(X + Z \leq t, Z > t)$$
$$+ P(X + Z > t, Z \leq t).$$

First, we will show that, as $t \nearrow \infty$,
$$P(X + Z > t, Z \leq t) \sim P(X > t)EX^-/(-EX).$$

Let $y_0 = \inf\{t \in \mathbb{R} : P(Z > t) < 1\}$. Then,
$$P(X + Z > t, Z \leq t)$$
$$= \frac{-1}{EX} \int_{y_0}^{t} P(X > t - s)P(X > s)\,ds + P(X > t - y_0)P(Z = y_0).$$

We now analyze the integral on the right-hand side of the previous display:
$$\int_{y_0}^{t} P(X > t-s)P(X > s)\,ds$$
$$= \int_{0}^{t-y_0} P(X > t - y_0 - s)P(X > s + y_0)\,ds$$
$$= \int_{0}^{t-y_0} P(X > t - y_0 - s)P(X > s)\,ds$$
$$+ \int_{0}^{t-y_0} P(X > t - y_0 - s)[P(X > s + y_0) - P(X > s)]\,ds.$$

Let us define by $I_1$ and $I_2$ the two last integrals on the right-hand side of the display above. Then, assumption A yields
$$I_1 = \int_{0}^{t-y_0} P(X > t - y_0 - s)P(X > s)\,ds$$
$$\sim 2P(X > t)EX^+ \qquad \text{as } t \nearrow \infty.$$

Now, for the integral $I_2$, we have

$$I_2 = \int_0^{t-y_0} P(X > t - y_0 - s) \, d\int_s^{s+y_0} P(X > u) \, du$$

$$= -\int_0^{t-y_0} \int_{t-y_0-s}^{t-s} P(X > u) \, du \, P(X \in ds)$$

$$+ P(X > 0) \int_{t-y_0}^t P(X > u) \, du - P(X > t - y_0) \int_0^{y_0} P(X > u) \, du.$$

Note that

$$P(X > t - s)y_0 \leq \int_{t-y_0-s}^{t-s} P(X > u) \, du$$

$$= t \int_1^{1+y_0/t} P(X > ut - s - y_0) \, du$$

$$\leq P(X > t - s - y_0)y_0.$$

Hence, by virtue of Assumption A, we have that, as $t \nearrow \infty$,

$$\int_0^{t-y_0} \int_{t-y_0-s}^{t-s} P(X > u) \, du \, P(X \in ds) \sim P(X > t)y_0 P(X > 0).$$

Similarly, we obtain that

$$P(X > 0) \int_{t-y_0}^t P(X > u) \, du \sim P(X > t)y_0 P(X > 0)$$

as $t \nearrow \infty$, which yields

$$I_2 \sim -P(X > t) \int_0^{y_0} P(X > s) \, ds.$$

Combining these estimates, we obtain

$$P(X + Z > t, Z \leq t)$$
$$\sim (I_1 + I_2)/(-EX) + P(X > t - y_0)P(Z = y_0)$$
$$\sim 2P(X > t)EX^+/(-EX) - P(X > t) \int_0^{y_0} P(X > s) \, ds/(-EX)$$
$$+ P(X > t - y_0)P(Z = y_0).$$

Since

$$P(Z = y_0) = 1 - \frac{1}{(-EX)} \int_{y_0}^{\infty} P(X > s) \, ds,$$

we have

$$P(X + Z > t; Z \leq t)$$
$$\sim P(X > t)[2EX^+ + (-EX) - EX^+]/(-EX)$$
$$= P(X > t)EX^-/(-EX).$$

On the other hand,

$$P(X + Z \leq t, Z > t) = -\frac{1}{EX} \int_t^\infty P(X \leq t - s)P(X > s) \, ds$$
$$= -\frac{1}{EX} \int_{-\infty}^0 P(X \leq s)P(X > t - s) \, ds$$
$$\sim P(X > t)EX^-/(-EX)$$

as $t \nearrow \infty$. This yields the proof of the result. $\square$

The constant $a_*$ that characterizes the region where inequality (11) holds will play an important role in the construction of our algorithm. The bound (11) indicates that on the region $(-\infty, a_*]$ the approximation to the zero-variance change-of-measure based on $v(\cdot)$ is good enough to control the variance of the likelihood ratio in our simulations. Finding $a_*$ can be done numerically or analytically depending on the problem at hand. For implementation, the simulator can choose any value of $\gamma$ (for instance, $\gamma = 1/2$) or optimize the asymptotic upper bound that we shall obtain in Theorem 3, which we now are ready to state and prove.

Consider the importance sampling change-of-measure generated by

$$
\begin{aligned}
Q_{a_*}(y, dz) &= \frac{P(y + X \in z + dz)v(z + a_*)}{w(y + a_*)} \\
&= P(y + X \in z + dz | Z + X \geq -y - a_*).
\end{aligned}
\tag{12}
$$

Then, we will show that the corresponding estimator defined as

$$R = I(\tau(0) < \infty) \prod_{j=1}^{\tau(0)} \frac{w(S_{k-1} + a_*)}{v(S_k + a_*)} \tag{13}$$

has bounded relative variance as $S_0 = y \searrow -\infty$.

THEOREM 3. *Suppose that Assumption* A *is in force. Fix* $\gamma \in (0, 1)$ *and select* $a_* = a_*(\gamma) \in (-\infty, 0]$ *as in* (11). *Then,*

$$E_y^{Q_{a_*}} R^2 \leq (1 - \gamma)^{-1} \kappa(a_*)^{-2} v(y + a_*)^2,$$

*where* $\kappa(a_*) = \inf_{z \geq 0}[v(z + a_*)] = P(Z > -a_*)$. *Consequently,*

$$\sup_{b > 0} E_{-b}^{Q_{a_*}}[R(b)^2]/P(M > b)^2 < \infty.$$

PROOF. Define
$$h(y) = I(y + a_* \leq 0) + (1 - \gamma)I(y + a_* > 0).$$

We wish to apply Proposition 2 so we must satisfy bound (7), which in our case can be written as

$$(14) \quad w(y + a_*)^{-1} Ev(X + y + a_*)h(X + y) \leq \left(\frac{v(y + a_*)}{w(y + a_*)}\right)^2,$$

for all $y \leq 0$. Here we have used the fact that $h(y) = 1$ for $y \leq 0$. Using the interpretation of $v(\cdot)$ as a tail probability, we note that the bound (14) can be expressed, for all $y \leq 0$, as

$$E(h(X + y) - 1 | X + Z > -y - a_*) \leq \frac{v(y + a_*)^2 - w(y + a_*)^2}{w(y + a_*)^2}.$$

Observe that
$$h(X + y) - 1 = -\gamma I(X \geq -y - a_*).$$

Therefore, it suffices to verify that, for all $y \leq 0$,

$$-\gamma P(X > -y - a_* | X + Z \geq -y - a_*)$$
$$\leq \frac{v(y + a_*)^2 - w(y + a_*)^2}{w(y + a_*)^2}.$$

However, it follows since $Z \geq 0$ and using the fact that $w(y) = P(X + Z \geq -y)$, that the previous inequality holds if and only if, for all $y \leq 0$,

$$-\gamma \leq \frac{v(y + a_*)^2 - w(y + a_*)^2}{P(X > -y - a_*)w(y + a_*)},$$

which is true by definition of $a_*$. The conclusion of the result follows directly from Propositions 2 and 3, the fact that $P(M > b) \sim v(-b + a_*)$ as $b \nearrow \infty$ and that the ratio $P(M > b)/v(-b + a_*)$ is bounded as a function of $b$ on compact sets. □

Our approach to the study of the issue of asymptotically vanishing relative error will involve taking advantage of extreme value theory; see, for instance, Embrechts, Klüppelberg and Mikosch (1997), Section 3.3. We say that $X_1$ belongs to the domain of attraction of $H$ [denoted by $X_1 \in MDA(H)$] if $H$ is nondegenerate and there exists a sequence of constants $c_n \geq 0$ and $d_n \in \mathbb{R}$ (for $n \geq 1$) such that

$$c_n^{-1}(\max(X_1, \ldots, X_n) - d_n) \Longrightarrow H$$

as $n \nearrow \infty$. The random variable $H$ must follow a so-called extreme value distribution which, due to the Fisher–Tippett theorem [see Embrechts, Klüppelberg and Mikosch (1997), page 121], can be of only three types. Only the cases when $H$ has Frechet distribution, given by

$$\Phi_\alpha(x) = \exp(-x^{-\alpha})I(x>0), \qquad \alpha > 0,$$

or when $H$ follows a Gumbel distribution described via

$$\Lambda(x) = \exp(-\exp(-x))$$

are of interest to us. The class $MDA(\Phi_\alpha)$ is precisely the class of regularly varying distributions with index $\alpha > 0$ [i.e., $P(X>x) = x^{-\alpha}L(x)$, where $L(\cdot)$ is slowly varying at infinity], whereas $MDA(\Lambda)$ contains other commonly used heavy-tailed distributions, such as log-normal and Weibull. The normalization constants in the definition of $H$ (i.e., the $c_n$'s and $d_n$'s) depend on the so-called auxiliary function, which is defined via

$$\xi(x) = \frac{\int_x^\infty P(X_1 > t)\,dt}{P(X_1 > x)}.$$

The following result of Asmussen and Kluppelberg (1996) provides some asymptotic properties of the zero-variance change-of-measure as $b \nearrow \infty$. These properties will be useful in verifying that our estimator possesses asymptotically vanishing relative variance as $b \nearrow \infty$.

THEOREM 4 [Asmussen and Kluppelberg (1996)]. *Assume either that $X_1$ is regularly varying with index $\alpha > 1$ or that assumption A holds and $X_1 \in MDA(\Lambda)$. Then, given $S_0 = -b < 0$ and $\tau(0) < \infty$,*

$$\left( \frac{S_{\tau(0)-1}}{\xi(b)}, \frac{\tau(0)}{\xi(b)}, \frac{S_{\tau(0)}}{\xi(b)} \right) \Longrightarrow (V, V/|EX|, T)$$

*as $b \nearrow \infty$, where $V$ and $T$ are a pair of random variables with joint distribution $P(V > x, T > y) = P(H > x + y)$.*

With the previous result in hand, we now can sharpen the conclusion of Theorem 3 to obtain asymptotically vanishing relative error. The next result provides theoretical justification for the empirical performance discovered in the numerical experiments shown in Section 5, in which the accuracy for a given number of runs is seen to improve when $b$ gets larger.

THEOREM 5. *Assume either that $X_1$ is regularly varying with index $\alpha > 1$ or that assumption A holds and $X_1 \in MDA(\Lambda)$. Then, one can choose $\gamma(b) \searrow 0$ and $-a_*(b) \nearrow \infty$ as $b \nearrow \infty$ such that the estimator provided by Algorithm 1 [defined as $R$ in (13)] satisfies*

$$\lim_{b \longrightarrow \infty} \frac{E_{-b}^{Q^*} R^2}{P(M>b)^2} = 1.$$

PROOF. Under the stated assumptions, $\xi(b) \nearrow \infty$ as $b \nearrow \infty$ [see Asmussen and Kluppelberg (1996)]. Furthermore, because both $X_1$ and $Z$ are long tailed, it follows that, for any constant $c > 0$, $\xi(b+c) \sim \xi(b)$ as $b \nearrow \infty$. We start by noting that

$$E_{-b}^{Q^*} R^2 = E_{-b}\left(I(\tau(0) < \infty) \prod_{j=1}^{\tau(0)} \frac{w(S_{k-1} + a_*)}{v(S_k + a_*)}\right)$$

$$= P(M > b)w(-b + a_*)$$

$$\times E_{-b}\left(\prod_{j=1}^{\tau(0)-1} \frac{w(S_k + a_*)}{v(S_k + a_*)} \frac{1}{v(S_{\tau(0)} + a_*)} \bigg| \tau(0) < \infty\right).$$

The Cauchy–Schwarz inequality implies that

$$E_{-b}\left(\prod_{j=1}^{\tau(0)-1} \frac{w(S_k + a_*)}{v(S_k + a_*)} \frac{1}{v(S_{\tau(0)} + a_*)} \bigg| \tau(0) < \infty\right)$$

(15)
$$\leq E_{-b}\left(\prod_{j=1}^{\tau(0)-1} \frac{w(S_k + a_*)^2}{v(S_k + a_*)^2} \bigg| \tau(0) < \infty\right)^{1/2}$$

$$\times E_{-b}\left(\frac{1}{v(S_{\tau(0)} + a_*)^2} \bigg| \tau(0) < \infty\right)^{1/2}.$$

Consider the first term in (15), which involves the ratios $w(S_k + a_*)^2 / v(S_k + a_*)^2$. We can again use a Lyapunov argument as the one introduced in the proof of Theorem 3. In fact, a completely analogous argument as the one given there shows that, for each $\gamma \in (0, 1)$, there exists a value of $a_* < 0$ for which

$$E_{-y}\left(\prod_{j=0}^{\tau(0)-1} \frac{w(S_k + a_*)^2}{v(S_k + a_*)^2} \bigg| \tau(0) < \infty\right) \leq \frac{1}{1 - \gamma}.$$

In fact, one just chooses $a_* < 0$ for which the inequality

(16)
$$-\gamma \leq \frac{v(y + a_*)^2 - w(y + a_*)^2}{P(X > -y - a_*)w(y + a_*)} \frac{u^*(y)P(X > -y - a_*)}{w(y + a_*)P(X > -y)}$$

$$= \frac{v(y + a_*)^2 - w(y + a_*)^2}{P(X > -y - a_*)w(y + a_*)} \frac{u^*(y)\xi(-y - a_*)v(y + a_*)}{v(y)\xi(-y)w(y + a_*)}$$

holds uniformly in $y \leq 0$. In view of Proposition 3, it suffices to analyze the ratio $\xi(-y)/\xi(-y - a_*)$. But because $\xi(b) \nearrow \infty$ as $b \nearrow \infty$ and $\xi(b+c) \sim \xi(b)$ as $b \nearrow \infty$, it follows that

$$\sup_{y \leq 0} \sup_{a \leq 0} \xi(-y)/\xi(-y - a) < \infty,$$

which is more than what is needed to guarantee that inequality (16) holds for $a_* < 0$ sufficiently negative.

With $a_*$ selected as above, observe that

$$
\begin{aligned}
E_{-b}&\left(\frac{1}{v(S_{\tau(0)} + a_*)^2}\bigg|\tau(0) < \infty\right) \\
&\leq P_{-b}(S_{\tau(0)} \geq -a_*|\tau(0) < \infty) + \frac{P_{-b}(S_{\tau(0)} \leq -a_*|\tau(0) < \infty)}{P(Z > -a_*)^2}.
\end{aligned}
\tag{17}
$$

By virtue of Theorem 4, we have that the right-hand side of (17) converges to 1 as $b \nearrow \infty$. We may therefore conclude that, given $\gamma \in (0,1)$, there exists a selection of $a_* > 0$ for which

$$\overline{\lim}_{b \to \infty} \frac{E_{-b}^{Q^*} R^2}{P(M > b)^2} \leq \frac{1}{1-\gamma}.$$

Since $\gamma > 0$ is arbitrary, we obtain the result by sending $\gamma \searrow 0$ and (possibly) also sending $a_*(\gamma) \searrow -\infty$ at a sufficiently slow rate. $\square$

REMARK. Although the previous result is intended only to provide a theoretical justification for the numerical performance found in our experiments, one can, in principle, find a computable constant $a_* < \infty$ for which

$$\overline{\lim}_{b \to \infty} \frac{E_{-b}^{Q^*} R^2}{P(M > b)^2} \leq \frac{1}{1-\gamma}.$$

This is clear from display (16) in the proof of Theorem 5. Note that everything in the left-hand side of (16) can, in principle, be evaluated, except of course $u^*(y)$. However, it suffices to find a computable bound for $u^*(y)/v(y)$, which can be obtained in many different ways, one of them through the use of the Lyapunov inequalities. Indeed, note that Theorem 3 and the fact that $E_{-b}^{Q^*} R^2 \geq u^*(-b)^2$ (which follows by the Cauchy–Schwarz inequality) could be used to obtain a computable upper bound for $\sup_{y \leq 0} u^*(y)/v(y)$.

**4. The algorithm and complexity analysis.** Recall that $Z$ is the nonnegative r.v. for which

$$P(Z > t) = \min\left[-(EX)^{-1} \int_t^\infty P(X > s)\,ds, 1\right].$$

The function $v(\cdot)$ is then defined through the relation $v(t) = P(Z > -t)$ and $w(\cdot)$ is given by $w(y) = P(X + Z > -y)$. We assume that $v$ and $w$ are either available in closed form or can be easily computed numerically. Note that in the conventional light-tailed setting, the sampler suggested through large deviations approximations to $v$ (and hence, to $w$) can often be implemented via "exponentially twisting" the increment distribution. Actual implementation of an importance sampler based on exponential twisting requires that the moment generating function be computable either in closed form or through a cheap numerical calculation, and that the correct twisting parameter (usually as characterized through a root of the gradient of the log moment generating function) be easily computable. Our assumptions on $v$ and $w$ can be viewed as the heavy-tailed analog to this requirement in the light-tailed case.

For fixed $\gamma \in (0,1)$, set $a_* = a_*(\gamma) \leq 0$ satisfying (11). We wish to estimate

$$u^*(-b) = P(\tau(0) < \infty | S_0 = -b),$$

for $b > 0$. Our proposed algorithm proceeds as follows.

Algorithm 1.

STEP 1. Initialize $s = -b$, $R = 1$.
STEP 2. Set $y \longleftarrow s$, generate a random variable $Y$ with law

$$P(Y \in t + dt) = P(X \in t + dt | X + Z > -y - a_*),$$

and update $s \longleftarrow y + Y$,

$$R \longleftarrow w(y + a_*)v(s + a_*)^{-1}R$$
$$= P(Z + X > -y - a_*)P(Z > -s - a_*)^{-1}R.$$

STEP 3. If $s > 0$ then return $R$ and STOP, otherwise, go to STEP 2.

Theorem 3 implies that the above algorithm is strongly efficient, in the sense that the number of simulation runs required to estimate $P(M > b)$ to a given relative accuracy is bounded in $b$. Within the rare-event simulation community, this statistical notion (and its close relative logarithmic efficiency) is the commonly accepted standard that a good algorithm should achieve.

However, a more demanding notion is to study the computational complexity of the algorithm. Roughly speaking, the goal is to show that the number of floating point operations required to compute $P(M > b)$ to a given relative accuracy increases at a slower rate than that associated with the use of crude Monte Carlo. [Note that it is typically unrealistic to expect that the number of floating point operations can be uniformly bounded over $b$, because the number of r.v.'s required to generate the random object $I(M > b)$ is increasing in $b$ in expectation]. In our current setting, the required number of floating point operations is determined by the number of simulation runs required to estimate $P(M > b)$ to a given relative accuracy multiplied by the expected number of floating point operations per run.

Since Theorem 3 has already established the boundedness of the number of simulation runs, the key issue then becomes estimating the expected number of floating point operations per run (as a function of $b$).

Note that a complete analysis of this issue is impossible without having a model for floating point arithmetic that attaches different costs to such computations as special function evaluations (e.g., numerically evaluating the exponential function), generating uniform random variates and performing comparison operations. In addition, such a complexity analysis requires explicit specification of the numerical effort involved in evaluating $v$ and $w$, both in the case in which closed forms are available and (even more critically) in the setting in which a numerical integration scheme is used to compute $w$. These issues arise even in the setting of light-tailed rare-event simulation, in which the algorithm typically depends on exponential twisting. In particular, the issue of numerical computation of the log moment generating function, its gradient and the associated roots would be a necessary element in a complete complexity analysis of a light-tailed rare-event simulation algorithm.

While a complete complexity analysis is no doubt a worthwhile undertaking, we present here a simplified analysis of what we believe reflect the key pragmatic complexity issues. We take the point of view that the expected number of floating point operations per run is roughly proportional to the total number of uniform random variables generated per run. The expected number of uniform random numbers required per run is obtained as the product of the number of steps needed by the importance sampler (having transition kernel $Q_{a_*}$) to cross level 0 from $-b$, multiplied by the typical number of uniform random variables required to generate each increment of the random walk as governed by the kernel $Q_{a_*}$. We shall argue, later in this section, that the expected number of steps needed by our importance sampler (having transition kernel $Q_{a_*}$) to cross level 0 from initial position $-b$ grows linearly in $b$; see Proposition 4. The question of how many uniform random variates are required, on average, per increment of $Q_{a_*}$ is very specific to the precise form of the distribution of $X$ and to the ingenuity employed in developing a variate generation scheme for simulating from $Q_{a_*}$'s increment distribution. To illustrate this point, we provide an acceptance-rejection algorithm later in this section that uses a bounded (in expectation) number of uniform random variates per increment simulated (that is bounded both in $b$ and the position $x$ of the sampler) whenever $X$ has a regularly varying continuous density (assuming that there exists a variate generation scheme that generates $X$ using a finite—in expectation—number of uniforms). It follows that the total number of uniform random variates required per simulation run grows at a linear rate in $b$. On the other hand, for crude Monte Carlo, the number of simulations runs required to compute $P(M > b)$ to a given relative accuracy scales in proportion to $1/P(M > b)$.

Because the paths on which $M < b$ take an infinite number of steps to generate, the expected number of floating point operations required per run is infinite. [If the simulation is terminated after $t$ steps with $t$ chosen so as to produce an estimator bias of a small (and given) magnitude relative to $P(M > b)$, both $t$ and the number of steps increase linearly in $b$.] Thus, our importance sampler provides a substantial improvement in computational complexity relative to crude Monte Carlo.

The following result (whose proof is given at the end of the section) provides sufficient conditions to ensure that Algorithm 1 terminates in at most $O(b)$ steps, given $S_0 = -b$.

PROPOSITION 4. *Assume that $E(X_1^p; X_1 > 0) < \infty$ for some $p > 1$ and that Assumption A is in place. Then*

$$E_{-b}^{Q_{a*}} \tau(0) = O(b)$$

*as $b \nearrow \infty$.*

In order to complete our complexity analysis, it is important to observe that the execution of STEP 2 of the algorithm involves a one dimensional rare-event type simulation problem. We have assumed that $v(\cdot)$ and $w(\cdot)$ can be easily evaluated. Nevertheless, it could be the case that finding explicitly the distribution of $Y$ in STEP 2 could be difficult or numerically expensive. We shall argue that the variates in STEP 2 can be simulated through a suitable acceptance/rejection scheme. Note, however, that one has to design the scheme in such a way that the acceptance probability remains uniformly bounded (in $y$) away from zero. By doing this, the generation of the random walk increments in STEP 2 under the importance sampling distribution has uniformly bounded complexity as $b \nearrow \infty$. Consequently, given Proposition 4, the expected number of variates required to run Algorithm 1 will be of order $O(b)$ as $b \nearrow \infty$.

Typically, acceptance/rejection schemes, such as those indicated in the previous paragraph, although not complicated, must be designed based on specific characteristics of the problem at hand. Assume that $X$ has a continuous density $f_X(\cdot)$. STEP 2 of Algorithm 1 requires sampling a r.v. $Y$ with density $f_Y(\cdot)$ defined, for $b \geq 0$, as

$$f_Y(z; b) = v(-b + z) f_X(z) / w(-b).$$

The objective is to find an easy way to simulate r.v. $\widetilde{Z}$ with computable density $f_{\widetilde{Z}}(z; b)$ such that, for all $z \in \mathbb{R}$,

(18) $$f_Y(z; b) \leq p_{acc}(b)^{-1} f_{\widetilde{Z}}(z; b),$$

where the acceptance probability, $p_{acc}(b)$, satisfies $\inf_{b \geq 0} p_{acc}(b) > 0$.

In order to illustrate the construction of $f_{\widetilde{Z}}(\cdot)$, let us assume that $f_X(\cdot)$ is regularly varying. We pick $\theta \in (0,1)$ and define

$$c(b) = P(X \leq b - \theta b)\frac{P(Z > \theta b)}{P(Z > b)} + \frac{P(X > b - \theta b)}{P(Z > b)},$$

$$\lambda_0(b) = c(b)^{-1} P(X \leq b - \theta b) P(Z > \theta b)/P(Z > b),$$

$$\lambda_1(b) = c(b)^{-1} P(X > b - \theta b)/P(Z > b).$$

Then, the mixture density

$$f_{\widetilde{Z}}(z; b) = \lambda_0(b) \frac{f_X(z) I(z \leq b - \theta b)}{P(X \leq b - \theta b)}$$

$$+ \lambda_1(b) \frac{f_X(z)}{P(X > b - \theta b)} I(z > b - \theta b)$$

satisfies

$$f_Y(z; b) \leq mc(b) f_{\widetilde{Z}}(z; b),$$

where

$$m \geq \sup_{b \geq 0}[P(Z > b)/P(Z + X > b)].$$

The acceptance probability using $f_{\widetilde{Z}}(z; b)$ as proposal is $[mc(b)]^{-1}$. Using elementary properties of regularly varying functions, it follows that $\inf_{b \geq 0}[c(b) \times m]^{-1} > 0$.

We conclude the section with a proof of Proposition 4.

LEMMA 1. *Suppose that Assumption A is in force and that $E(X_1^p; X_1 > 0) < \infty$ for some $p > 1$. Then there exists $t_0 > 0$ and $\varepsilon > 0$ such that, for all $t \geq t_0$,*

$$E[X|X + Z > t] \geq \varepsilon.$$

PROOF. The assumptions imply that $X_1$ and $Z$ must be subexponential. In particular, it follows that $P(X + Z > t) \sim P(Z > t)$ as $t \nearrow \infty$. Thus, it suffices to show that

(19) $$\varliminf_{t \to \infty} \left( \frac{-E[X^-; X + Z > t]}{P(Z > t)} + \frac{E[X^+; X + Z > t]}{P(Z > t)} \right) > 0.$$

The Bounded Convergence Theorem implies that

$$\frac{-E[X^-; X + Z > t]}{P(Z > t)} = \int_{-\infty}^0 s \frac{P(Z > t - s)}{P(Z > t)} P(X \in ds) \longrightarrow -EX^-$$

as $t \nearrow \infty$. On the other hand, we have that

$$E[X^+; X + Z > t]$$
$$= E[X; X + Z > t; X \geq 0]$$
$$= E \int_0^\infty I(X + Z > t; X \geq 0; X \geq s) \, ds$$
$$= \int_0^\infty P(X + Z > t; X \geq s) \, ds$$
$$= \int_0^{t-y_0} P(X + Z > t; X \geq s) \, ds + P(Z > t - y_0)(-EX),$$

where $y_0 = \inf\{t \in \mathbb{R} : P(Z > t) < 1\}$. Now,

$$\int_0^{t-y_0} P(X + Z > t; X \geq s) \, ds$$
$$= \int_0^{t-y_0} P(X + Z > t; s \leq X \leq t - y_0) \, ds$$
$$\quad + (t - y_0) P(X > t - y_0).$$

The first integral in the right-hand side of the previous display is greater or equal to

$$\int_0^{t-y_0} P(Z > t - s) P(s \leq X \leq t - y_0) \, ds \sim P(Z > t) EX^+.$$

On the other hand, it follows that

$$tP(X > t) \geq \int_t^{2t} P(X > s) \, ds$$
$$= [P(Z > t) - P(Z > 2t)](-EX).$$

Observe that if $E(X_1^p; X_1 > 0) < \infty$ for some $p > 0$, then there exists $\delta > 0$ such that the map $t \longrightarrow t^\delta P(Z > t)$ is eventually decreasing. Therefore,

(20) $\qquad \overline{\lim}_{t \to \infty} P(Z > 2t)/P(Z > t) \leq (1/2)^\delta < 1.$

Putting all the previous estimates together (and using the fact that $Z$ has a long tail), we obtain that the limit in (19) is greater or equal to

$$-EX^- - EX + EX^+ - (1 - (1/2)^\delta) EX$$
$$= -(1 - (1/2)^\delta) EX > 0,$$

which is more than we need in order conclude the proof of the lemma. □

Finally, we provide the proof of Proposition 4.

PROOF OF PROPOSITION 4. It follows from Lemma 1 and Chebyshev's inequality that there exists $a < 0$ and $\varepsilon > 0$ such that

$$\sup_{y \leq a} E(X | X + Z > -y - a) > \varepsilon. \tag{21}$$

Now, set $\tau(a) = \inf\{n \geq 1 : S_n > a\}$. It follows from (21) then that on $\{\tau(a) > n\}$, there exists $\varepsilon > 0$ such that

$$E^{Q_{a*}}(S_{n+1} | S_n) - S_n > \varepsilon$$

and, therefore [letting $\min(n, \tau(0)) \triangleq n \wedge \tau(0)$], it is not hard to see that $M_n = S_{n \wedge \tau(0)} - (n \wedge \tau(a))\varepsilon$ is a submartingale (under $Q_{a_*}$). In particular, we obtain that

$$\varepsilon E_y^{Q_{a*}}[n \wedge \tau(a)] \leq E_y^{Q_{a*}} S_{n \wedge \tau(a)} - y \leq -y.$$

Finally, the monotone convergence theorem yields

$$E_y^{Q_{a*}} \tau(a) \leq |y|/\varepsilon.$$

On the other hand, we have that

$$\sup_{y \in [a,0]} E_y^{Q_{a*}} \tau(a) \leq 1 - \frac{1}{\varepsilon} \sup_{a \leq y \leq 0} E[X + y; X + y \leq a | X + Z > -y - a_*] < \infty.$$

Therefore, it follows from a geometric trials argument that

$$E_{-b}^{Q_{a*}} \tau(0)$$
$$\leq E_{-b}^{Q_{a*}} \tau(a) + \left[\sup_{a \leq y \leq 0} P(X > -a | X + Z > -y - a_*)\right]^{-1} \sup_{y \in [a,0]} E_y^{Q_{a*}} \tau(a)$$
$$\leq b \cdot m$$

for some $m \in (0, \infty)$, which yields the proof of the result. □

**5. Empirical validation.** We first consider a class of models for which other provably efficient algorithms have been developed. This permits a direct computational comparison of our algorithm against other existing methods. In particular, we shall consider here an $M/G/1$ queue with traffic intensity $\rho = 1/2$ and with Pareto service times having tail $P(V > t) = (1 + t)^{-2.5}$ (we use $V$ to denote a generic service time). As noted in the Introduction, two existing competing algorithms for this class of models are the ones proposed by Asmussen and Kroese (2006) (AK) and Dupuis, Leder and Wang (2006) (DLW). AK's procedure was designed to deal with regularly varying and Weibull-type service times, whereas DLW's works only for regularly varying distributions. There is only one other algorithm that can be shown to be strongly efficient for a single-server $G/G/1$ queue with arbitrary renewal arrivals, due to Blanchet, Glynn and Liu (2007) (BGL). It should be noted that it requires regularly varying service times, and hence, does not cover the range of $G/G/1$ models covered by this paper's algorithm. The BGL procedure does not exploit the representation of the steady-state distribution of the $G/G/1$ queue in terms of the maximum of a random walk, but instead takes advantage of the regenerative ratio representation for the steady-state distribution.

We have also added a computational comparison against a more refined implementation of the algorithm introduced in this paper, which takes advantage of the exponential inter-arrival times to help speed up the time it takes for the random walk to hit the target set. This refined implementation can be found in Blanchet and Li (2006) (BL), and relies on the fact that, for an $M/G/1$ queue, the distribution of the ladder heights can be computed explicitly. We therefore can apply the algorithm of the current paper to a first-passage time problem involving a strictly increasing random walk that is killed at a geometric time, thereby basically running the algorithm at the level of ladder epochs (and saving the computer time associated with generating the transitions between ladder epochs that occurs in the algorithm described in this paper).

Table 1, based on 10,000 samples, compares the performance of our algorithm (which is denoted in the table below by BG) against the AK, DLW and BGL procedures. In our algorithm, we chose $a_* = 10$ and carefully implemented a numerical integration routine in order to compute $w(\cdot)$; $v(\cdot)$ can be evaluated in closed form in terms of an incomplete gamma function. The sampling of each of the increments was done using an acceptance rejection procedure similar to the one explained in Section 4 right after display (18).

Perhaps not surprisingly, given that a closely related variant of our estimator (in which $a_*$ can increase with $b$) exhibits asymptotically vanishing relative error (as indicated by Theorem 5), one can see that the coefficient of variation displayed for BG and BL diminishes as the level increases. As indicated above, the advantage of the BL implementation over the algorithm discussed here (BG) is that BG requires $O(b)$ variate generations per replication of the estimator, whereas the requirement is $O(1)$ as $b \nearrow \infty$ for BL—of course, such speed up relies heavily on the assumption of Poisson arrivals. It should be noted that the AK and DLW algorithms also require $O(1)$ variates per replication for their estimators, and enjoy relatively simple implementations, particularly for the case of the AK estimator. Finally, we note that the BGL procedure also takes $O(b)$ variate generations per replication of the estimator.

We conclude this section with a problem instance for which BG is the only currently available procedure for efficient tail estimation. In particular, consider a $G/G/1$ queue with deterministic inter-arrival times equal to 1 and a service time tail distribution given by

$$P(V > t) = \exp(-2t^{1/2}).$$

TABLE 1
*Numerical estimates for $u(S_0)$ with Pareto tails and $\rho = 1/2$*

| [Estimation]<br>[Std. error]<br>[Conf. interval] | $S_0 = -10^2$ | $S_0 = -10^3$ | $S_0 = -10^4$ |
|---|---|---|---|
| $v(S_0)$ | $9.663e-03$ | $3.151e-05$ | $9.996e-07$ |
| AK | $5.995e-04$ | $3.145e-05$ | $9.980e-07$ |
|  | $7.395e-06$ | $2.186e-07$ | $6.945e-09$ |
|  | $[5.850e-04,$ | $[3.102e-05,$ | $[9.844e-07,$ |
|  | $6.140e-04]$ | $3.188e-05]$ | $1.012e-06]$ |
| BG | $5.485e-04$ | $3.146e-05$ | $9.980e-07$ |
|  | $2.984e-06$ | $9.725e-08$ | $2.073e-09$ |
|  | $[5.427e-04,$ | $[3.126e-05,$ | $[9.939e-07,$ |
|  | $5.543e-04]$ | $3.165e-05]$ | $1.002e-06]$ |
| BL | $9.750e-04$ | $3.162e-05$ | $9.980e-07$ |
|  | $4.363e-06$ | $1.982e-07$ | $4.511e-09$ |
|  | $[9.664e-04,$ | $[3.123e-05,$ | $[9.892e-07,$ |
|  | $9.836e-04]$ | $3.201e-05]$ | $1.007e-06]$ |
| BGL | $1.022e-03$ | $3.167e-05$ | $1.128e-06$ |
|  | $3.835e-05$ | $1.598e-06$ | $7.280e-08$ |
|  | $[9.468e-04,$ | $[2.854e-05,$ | $[9.853e-07,$ |
|  | $1.097e-03]$ | $3.480e-05]$ | $1.271e-06]$ |
| DLW | $1.046e-03$ | $3.163e-05$ | $9.905e-07$ |
|  | $5.195e-06$ | $1.694e-07$ | $2.993e-09$ |
|  | $[1.036e-03,$ | $[3.130e-05,$ | $[9.846e-07,$ |
|  | $1.056e-03]$ | $3.196e-05]$ | $9.964e-07]$ |

The increment $X_n$ has a distribution given by $X = V - 1$ so that

$$EX = \int_0^\infty \exp(-2t^{1/2})\,dt - 1$$
$$= 1/2 - 1 = -1/2,$$

which implies that the traffic intensity is $\rho = 1/2$.

To run the algorithm, we selected $a_* = -10$. We also tried $a_* = -5$ and $-55$. For $a_* = -5$, the sample coefficient of variation was slightly lower than the one that we display below, but not too much. For $a_* = -55$, we obtained sample coefficients of variation no larger than 100. In both cases, the corresponding pointwise estimates were very consistent with those displayed below.

The most interesting part of the implementation involves Step 2, namely, sampling from the r.v. $Y$ with law

$$P(Y \in t + dt) = P(X \in t + dt | X + Z > \beta).$$

For this step we use an acceptance/rejection procedure. Again, we make sure that the acceptance probability remains uniformly bounded away from zero over as $\beta \nearrow \infty$. Let $m = \lfloor \beta^{1/2} \rfloor$ and note that

$$P(X \in t + dt | X + Z > \beta)$$
$$= \frac{1}{w(-\beta)} f_X(t) P(Z > \beta - t)$$
$$\leq f_X(t) I(-1 \leq t < 0) \frac{P(Z > \beta)}{w(-\beta)}$$
$$+ \sum_{k=0}^{m-1} f_X(t) I(k\beta^{1/2} \leq t \leq (k+1)\beta^{1/2}) \frac{P(Z > \beta - (k+1)b^{1/2})}{w(-\beta)}$$
$$+ f_X(t) I(m\beta^{1/2} \leq t \leq \beta) \frac{P(Z > 0)}{w(-\beta)} + f_X(t) I(t \geq \beta) \frac{1}{w(-\beta)}.$$

Using the dominating density induced by the expression in the right-hand side, we have that Step 2 can be performed in $O(\beta^{1/2})$ operations [hence, a single sample path generated by the proposed algorithm takes at most $O(b^{3/2})$ operations].

Table 2 illustrates the performance of the algorithm. BG is the estimator based on our importance sampling scheme using 20,000 replications. In order to validate the implementation of the algorithm, we constructed, for $S_0 = -10$, a crude Monte Carlo estimator based on 500,000 replications. The estimator was obtained using the regenerative ratio formula [see Asmussen (2003), page 268, equation (1.6)]. An approximate 95% confidence interval for $u(-10)$ based on these 500,000 samples is $[1.862e - 02, 2.179e - 02]$ (the point estimate was $2.020e - 02$). It is worth noting that the width of our importance sampling confidence interval is about 1/2 of that produced by crude Monte Carlo, with 25 times fewer samples [for a probability that is just moderately small, as is the case of $u(-10)$]. We did not apply crude Monte Carlo at the other values of $S_0$, because of the prohibitive amount of computation required.

The column CV reports the estimated coefficient of variation of our estimator, that is, the (estimated) standard deviation divided by the sample mean.

TABLE 2
*Numerical estimates for $u(S_0)$ with Weibull tails and $\rho = 1/2$*

| $S_0$ | $v(S_0)$ | BG | CV | 95% Conf. interval |
|---|---|---|---|---|
| $-10$ | $1.004e - 02$ | $1.942e - 02$ | 3.68 | $[1.857e - 02, 2.027e - 02]$ |
| $-50$ | $9.577e - 06$ | $1.783e - 05$ | 2.40 | $[1.724e - 05, 1.842e - 05]$ |
| $-250$ | $5.666e - 13$ | $7.076e - 13$ | 2.39 | $[6.842e - 13, 7.310e - 13]$ |
| $-500$ | $1.655e - 18$ | $1.897e - 18$ | 3.79 | $[1.797e - 18, 1.997e - 18]$ |
| $-650$ | $3.584e - 21$ | $3.971e - 21$ | 2.83 | $[3.815e - 21, 4.127e - 21]$ |

**Acknowledgments.** Thanks to Chenxin Li for his help with the numerical implementation of the routines, Bert Zwart for pointing out the need for condition $S^*$ in our proof of Proposition 3 and to Ton Dieker for his valuable comments on an earlier version of this paper. We are also grateful for the suggestions of the referee which helped improve the presentation of this paper.

## REFERENCES

ADLER, J., FELDMAN, R. and TAQQU, M., eds. (1998). *A Practical Guide to Heavy Tails*: *Statistical Techniques and Applications*. Birkhäuser, Boston. MR1652283

ASMUSSEN, S. (2003). *Applied Probability and Queues*. Springer, New York. MR1978607

ASMUSSEN, S. and GLYNN, P. (2007). *Stochastic Simulation*: *Algorithms and Analysis*. Springer, New York. MR2331321

ASMUSSEN, S. and BINSWANGER, K. (1997). Simulation of ruin probabilities for subexponential claims. *Ast. Bulletin* **27** 297–318.

ASMUSSEN, S., BINSWANGER, K. and HOJGAARD, B. (2000). Rare event simulation for heavy-tailed distributions. *Bernoulli* **6** 303–322. MR1748723

ASMUSSEN, S. and KLUPPELBERG, C. (1996). Large deviation results for subexponential tails, with applications to insurance risk. *Stoch. Proc. Appl.* **64** 103–125. MR1419495

ASMUSSEN, S. and KROESE, D. (2006). Improved algorithms for rare event simulation with heavy tails. *Adv. Appl. Probab.* **38** 545–558. MR2264957

BASSAMBOO, A., JUNEJA, S. and ZEEVI, A. (2006). On the efficiency loss of state-independent importance sampling in the presence of heavy-tails. *Oper. Res. Lett.* **34** 521–531.

BLANCHET, J., GLYNN, P. and LIU, J. C. (2007). Fluid heuristics, Lyapunov bounds and efficient importance sampling for a heavy-tailed G/G/1 queue. *QUESTA* **57** 99–113.

BLANCHET, J. and LI, C. (2006). Efficient rare-event simulation for geometric sums. *Proc. RESIM, Bamberg.* Germany.

BOROVKOV, A. A. and BOROVKOV, K. A. (2001). On probabilities of large deviations for random walks. I: Regularly varying distribution tails. *Theory Probab. Appl.* **49** 189–205.

BUCKLEW, J. (2004). *Introduction to Rare-Event Simulation*. Springer, New York. MR2045385

DUPUIS, P. and WANG, H. (2004). Importance sampling, large deviations, and differential games. *Stochastics Stochastics Rep.* **76** 481–508. MR2100018

DUPUIS, P., LEDER, K. and WANG, H. (2006). Importance sampling for sums of random variables with regularly varying tails. *TOMACS* **17**.

EMBRECHTS, P., KLÜPPELBERG, C. and MIKOSCH, T. (1997). *Modelling Extremal Events for Insurance and Finance*. Springer, New York. MR1458613

HAMMERSLEY, J. and MORTON, K. (1954). Poor man's Monte Carlo. *J. Roy. Statist. Soc. Ser. B* **16** 23–38. MR0064475

JUNEJA, S. and SHAHABUDDIN, P. (2002). Simulating heavy-tailed processes using delayed hazard rate twisting. *ACM TOMACS* **12** 94–118.

JUNEJA, S. and SHAHABUDDIN, P. (2006). Rare event simulation techniques: An introduction and recent advances. In *Handbook on Simulation* (S. Henderson and B. Nelson, eds.) 291–350. North-Holland, Amsterdam.

LIU, J. (2001). *Monte Carlo Strategies in Scientific Computing*. Springer, New York. MR1842342

MEYN, S. and TWEEDIE, R. (1993). *Markov Chains and Stochastic Stability*. Available at http://decision.csl.uiuc.edu/~meyn/pages/book.html. MR1287609

Rosenbluth, M. and Rosenbluth, A. (1955). Monte Carlo calculation of the average extension of molecular chains. *J. Chem. Phys.* **23** 356–359.

Siegmund, D. (1976). Importance sampling in the Monte Carlo study of sequential tests. *Ann. Statist.* **4** 673–684. MR0418369

Department of Statistics  
Harvard University  
Cambridge, Massachusetts 02138  
USA  
E-mail: blanchet@stat.harvard.edu

Department of Management  
Science and Engineering  
Stanford University  
Stanford, California 94305  
USA  
E-mail: glynn@stanford.edu